\documentclass{article}
\pdfoutput=1 
\font\myfont=cmr12 at 18pt
\newcommand{\TheTitle}{\myfont Delay compartment models from a stochastic process}

\usepackage{graphicx}
\usepackage{amsmath}
\DeclareMathOperator*{\argmin}{arg\,min}
\usepackage{amsfonts}
\usepackage{Styles/natbib}
\usepackage[toc,page]{appendix}
\graphicspath{{Figures/}} 
\usepackage{float}
\usepackage{xcolor}
\usepackage{enumerate}
\usepackage{chngcntr}
\usepackage{subcaption}
\title{{\TheTitle}}

\author{
  Christopher N. Angstmann\thanks{School of Mathematics and Statistics, UNSW Australia Sydney, NSW 2052, Australia} \and Anna V. McGann\footnotemark[1] \and Zhuang Xu\footnotemark[1] 
     }

\usepackage{amsopn}

\providecommand{\msc}[1]{\textbf{MSC:} #1}

\begin{document}
	
	\newcommand{\calL}{\mathcal{L}}

	\maketitle

	\begin{abstract}
		Compartment models with delay terms are widely used across a range of disciplines. The motivation to include delay terms varies across different contexts. In epidemiological and pharmacokinetic models, the delays are often used to represent an incubation period. In this work, we derive a compartment model with delay terms from an underlying non-Markov stochastic process. Delay terms arise when waiting times are drawn from a delay exponential distribution. This stochastic process approach allows us to preserve the physicality of the model, gaining understanding into the conditions under which delay terms can arise.  By providing the conditions under which the delay exponential function is a probability distribution, we establish a critical value for the delay terms. An exact stochastic simulation method is introduced for the generalized model, enabling us to utilize the simulation in scenarios where intrinsic stochasticity is significant, such as when the population size is small. We illustrate the applications of the model and validate our simulation algorithm on examples drawn from epidemiology and pharmacokinetics.\\		
		% \subclass{} %to be filled late
	\end{abstract}
\\
	\msc{33E30, 34K99, 60K15, 92C45, 92D30}
	
	\section{Introduction}
	Compartment models are frequently used for mathematical modelling across a broad range of disciplines, including pharmacokinetics \cite{T1937,N2003}, ecology \cite{M1994}, population dynamics \cite{H1975,M2000} and epidemic modelling \cite{KM1927,H2000}. These models split a phenomena into discrete states, treating each one as a compartment, and describe the dynamics between the compartments.
	
	Compartments models are represented by a set of coupled ordinary differential equations (ODEs). The populations of the compartments are considered, with terms describing interactions between the compartments, birth or removal processes. Members of each compartment are treated indistinguishably to simplify the model. This leads to every compartment being well-mixed and homogenous. 
	
	Two common compartment model representations are populations and concentrations. A population compartment model can be turned into a concentration compartment model by rescaling parameters \cite{HD1995}. Pharmacokinetic models tend to be formulated in terms of concentrations \cite{N2003}. Epidemiological models vary in their formulating, sometimes using populations or concentrations \cite{HD1995}.
	
	A limitation of ODE compartment models is that the dynamics are governed by their current state. However, there are many processes whose dynamics are dependent on a prior state, or even the entire history of the system. To represent such a history effect in a compartment model, it has become increasingly popular to include non-local operators in the governing equations \cite{AEHMMN2017,AEHMMN2021,S2011,RHRC2018}. These non-local operators lead to systems featuring fractional order derivatives or delay differential equations (DDEs). However, this has typically been done without consideration to preserving the underlying assumptions of the model such as: dimensional agreement, interpretation of the parameters and conservation of mass \cite{DM2009,DMM2010}. 
	
	In Angstmann et al. \cite{AEHMMN2017,AEHMMN2021} a general framework was established to incorporate fractional derivatives into compartment models while preserving the model's physicality. This was derived by considering the underlying stochastic process of individuals moving through compartments. By considering a particular waiting time in  a compartment, fractional derivatives naturally arise. 
	There has also been progress made to generalize this framework to produce a broader class of models \cite{ZWR2024}.  Similarly, our work seeks to establish a generalized framework for the incorporation of time-delays into compartment models as well as provide a stochastic simulation of compartment models with time-delays.
	
	Compartment models represented by DDEs span a wide range of disciplines \cite{RHRC2018}. These include ecology \cite{LS2000}, population dynamics \cite{H1948,AF1990}, epidemic modelling \cite{C1979,GAG2011} and predator-prey models \cite{AAdS1977}. The time-delay is typically motivated to represent an incubation time, such as accounting for infancy before an individual in the population can reproduce when modelling population growth \cite{TM1970} or allowing an infection to incubate before it can be transmitted to another individual \cite{HTMW2010}. A well established property of DDEs is that they often induce oscillatory solutions \cite{AGJ1984,G1986,K1993,G2013}. This has also served as a motivation to use a DDE model when the known behaviours of the system include oscillations \cite{GAG2011}.  
	
	Attention has been given to the technical aspects of DDE models. This has included solving initial value problems for DDEs \cite{KF2002,K1993,G2013}. As well as performing stability analyses \cite{C1979,HTMW2010}. However, insight into the modelling aspects by the consideration of an underlying stochastic process that results in a DDE.
	
	We derive a generalized compartment model from a continuous time random walk (CTRW) \cite{MW1965} in Section \ref{sec:deriv} and establish the conditions under which delay terms can arise. A derivation of a general compartment model with a non-Markovian removal process was established in Angstmann et al. \cite{AEHMMN2017,AEHMMN2021}.
	By taking this approach, we develop a further understanding into the dynamics that lead to time-delays, and ensure we preserve the physicality of the model. In Section \ref{sec:model_time_delay}, we show that a delay arises naturally if we select a delay exponential distribution as the non-Markovian waiting time for the removal process out of a compartment.
	
	A delay exponential function has previously been defined \cite{KIS2005,DK2006,ABHJX2023} along with many of its properties formalized \cite{ABHJX2023}. In Section \ref{sec:dexp_dist}, we introduce a delay exponential probability distribution from the delay exponential function. We show the required constraints on the delay exponential function under which it is a probability distribution.
	
	We also provide an exact simulation method for the resulting stochastic compartment model with delays in Section \ref{sec:simulation}, by making use of the general simulation framework established in Xu et al. \cite{XAHBHJ2023}. The stochastic simulation method was established for compartment models with fractional derivatives, and allows for a general non-Markovian removal process. Here we consider the case when the non-Markovian removal process is defined by a delay exponential distribution. In Section \ref{sec:example}, we apply the simulation method to pharmacokinetic and epidemiological models, and numerically demonstrate the convergence between stochastic and deterministic dynamics. We conclude with a discussion in Section \ref{sec:summary}.
	
	\section{Derivation of generalized compartment model} \label{sec:deriv}
	In order to develop a multi-compartment model with time delays, we will first consider a generalized compartment model. We consider a model with $N$ compartments, resulting in $N$ coupled differential equations. In general, the evolution of each compartment will be defined by the flux into the compartment and the removal processes out of the compartment at time, $t$. The specifics of each equation are system dependent.
	
	To begin the derivation, we consider an arbitrary compartment, $i$, within this model.
	The total flux into the $i$ compartment, $q_i(t)$, is the sum of $N_C$ distinct fluxes $\beta_{i,j}(t)$,
	\begin{equation}
		q_i(t)=\sum_{j=1}^{N_c} \beta_{i,j}(t).
	\end{equation}
	Our model will presume that the compartment is empty prior to $t=0$. Similarly, we consider the flux as starting with a burst at $t=0$, and then a continuous function for $t>0$. So, we can write the flux as,
	\begin{equation} \label{eq:flux}
		q_i(t)=i_0\delta(t) + q_i^+(t),
	\end{equation}
	where $i_0$ is the initial flux into the compartment and $q_i^+(t)$ is right continuous at $t=0$, and continuous for $t>0$. 
	Once a particle has entered this compartment through the flux it remains in the compartment until it is removed. We allow for $N_R$ Markovian removal processes. Each individual rate is written as $\lambda_{i,j}(t)$. The total Markovian removal rate can then be defined as,
	\begin{equation}\label{eq:omega_lambda}
		\omega_i(t)=\sum_{j=1}^{N_R}\lambda_{i,j}(t).
	\end{equation}
	The survival function for the Markovian removal processes can be given by,
	\begin{equation}\label{eq:Phi_def}
		\Phi_i(t,t_0)=\exp\left(-\int_{t_0}^t \omega_i(s)ds\right).
	\end{equation}
	As such it possesses the semi-group property,
	\begin{equation}
		\Phi_i(t,t_0) = \Phi_i(t,u)\Phi_i(u,t_0),
	\end{equation}
	where $t_0\leq u \leq t$. In addition, a non-Markovian removal process characterized by the survival function $\Psi_i(t-t_0)$ is included. To simplify the presentation, we consider only a single removal process of this form.
	
	Now, we can consider the expected number of particles in compartment $i$ at time $t$. For a particle to be in compartment $i$ at time $t$, it must have entered the compartment at an earlier time $t_0$ and survived all of the removal processes until time $t$, which we have taken to be independent. Hence, the expected number of particles in the compartment at time $t$ can be written as,
	\begin{equation}\label{eq:rho_t}
		\rho_i(t)=\int_0^t q_i(t_0)\Phi_i(t,t_0)\Psi_i(t-t_0)dt_0.
	\end{equation}
	We will consider the form of the non-Markovian removal process that leads to time-delay terms in the system in Section \ref{sec:model_time_delay}. 
	
	Differentiating Eq. \eqref{eq:rho_t}, we get the governing equation for an ensemble of a particles in the $i$th compartment \cite{AEHMMN2017,AEHMMN2021},
	\begin{equation} \label{eq:gov_single_compartment}
		\frac{d\rho_i}{dt}=q_i^+(t)-\omega_i(t)\rho_i(t)-\int_0^t \psi_i(t-t_0)\Phi_i(t,t_0)q_i(t_0)dt_0,
	\end{equation}
	where the waiting time density, $\psi_i(t)$, is related to the survival function by
	\begin{equation}
		\Psi_i(t) = 1-\int_{0}^{t}\psi_i(u)du.
	\end{equation}
	We will define the outgoing flux by the non-Markovian process, $F_{i}(t)$, as
	\begin{equation} \label{eq:F_def}
		F_i(t)=\int_0^t \psi_i(t-t_0)\Phi_i(t,t_0)q_i(t_0)dt_0.
	\end{equation}
	Using the semi-group property, Eqs. \eqref{eq:rho_t} and \eqref{eq:F_def} can be rewritten as,
	\begin{align}
		\frac{\rho_i(t)}{\Phi_i(t,0)}&=\int_0^t \Psi_i(t-t_0)\frac{q_i(t_0)}{\Phi_i(t_0,0)}dt_0,\\
		\frac{F_i(t)}{\Phi_i(t,0)}&=\int_0^t \psi_i(t-t_0)\frac{q_i(t_0)}{\Phi_i(t_0,0)}dt_0.
	\end{align}
	Then their Laplace transforms, with respect to $t$ can be written as,
	\begin{align}
		\mathcal{L}_t\left\{\frac{\rho_i(t)}{\Phi_i(t,0)}\right\}&=\mathcal{L}_t\left\{\Psi_i(t)\right\}\mathcal{L}_t\left\{\frac{q_i(t)}{\Phi_i(t,0)}\right\},\\
		\mathcal{L}_t\left\{\frac{F_i(t)}{\Phi_i(t,0)}\right\}&=\mathcal{L}_t\left\{\psi_i(t)\right\}\mathcal{L}_t\left\{\frac{q_i(t)}{\Phi_i(t,0)}\right\}.
	\end{align}
	These equations can be used to write,
	\begin{equation}
		\mathcal{L}_t\left\{\frac{F_i(t)}{\Phi_i(t,0)}\right\}=\frac{\mathcal{L}_t\left\{\psi_i(t)\right\}}{\mathcal{L}_t\left\{\Psi_i(t)\right\}}\mathcal{L}_t\left\{\frac{\rho_i(t)}{\Phi_i(t,0)}\right\}.
	\end{equation}
	Which substituted into Eq. \eqref{eq:gov_single_compartment} produces the governing equation with a general non-Markovian recovery rate is,
	\begin{equation}\label{eq:gov_comp_kern}
		\frac{d\rho_i}{dt}=q_i^+(t)-\omega_i(t)\rho_i(t)-\Phi_i(t,0)\int_0^t K_i(t-t_0)\frac{\rho_i(t_0)}{\Phi_i(t_0,0)}dt_0,
	\end{equation}
	where the memory kernel is defined as:
	\begin{equation} \label{eq:memory_kernel}
		\calL_t\{ K_i(t) \} = \frac{ \calL_t\{\psi_i(t)\} }{ \calL_t\{\Psi_i(t)\} } \,.
	\end{equation}
	Using different waiting times will lead to different governing equations. Taking $\Psi_i(t)$ to be an exponential survival function leads to a standard ODE model. In this case there is no non-Markovian survival function. Meanwhile, when $\Psi_i(t)$ is a Mittag-Leffler survival function, the model will feature fractional derivatives \cite{AEHMMN2017,AEHMMN2021}. We will now focus on the choice of survival function that leads to a DDE model. The necessary form of the survival function is a delay exponential distribution.                                
	
	\section{Delay exponential function as a survival function}\label{sec:dexp_dist}
	The delay exponential function is a delay analogue of the standard exponential function. It is defined through a truncated power series \cite{ABHJX2023} as,
	\begin{equation} \label{eq:dexp}
		\mathrm{dexp}(-\mu t;-\mu\tau)=\sum_{n=0}^\infty  \frac{(-\mu)^n(t-n\tau)^n}{\Gamma(n+1)}\Theta\left(\frac{t}{\tau}-n\right),\quad \frac{t}{\tau}\in  \mathbb{R}.
	\end{equation}
	Here $\mu>0$ is a scale parameter, $\tau\ge 0$ is the constant delay, and $\Theta(t)$ is the Heaviside function defined as, 
	\begin{equation} \label{eq_heaviside}
		\Theta(t)=\begin{cases} 0 & t<0,\\ 1& t\ge 0. \end{cases}
	\end{equation}
	Note that the power series representation of the exponential function is recovered when the limit as $\tau\rightarrow 0$ is taken.
	
	The choice of survival function that leads to delay terms in our compartment model is the delay exponential function:
	\begin{equation}\label{eq:delay_survival}
		\Psi(t)= \mathrm{dexp}(-\mu t;-\mu\tau),
	\end{equation}
	
	The corresponding density function of Eq. \eqref{eq:delay_survival}, $\psi(t)$, is:
	\begin{equation}\label{eq:dexpPDF}
		\psi(t)=\mu\mathrm{dexp}(-\mu(t-\tau);-\mu\tau).
	\end{equation} 
	
	\begin{figure}[H]
		\begin{center}
			\includegraphics[width=0.9\linewidth]{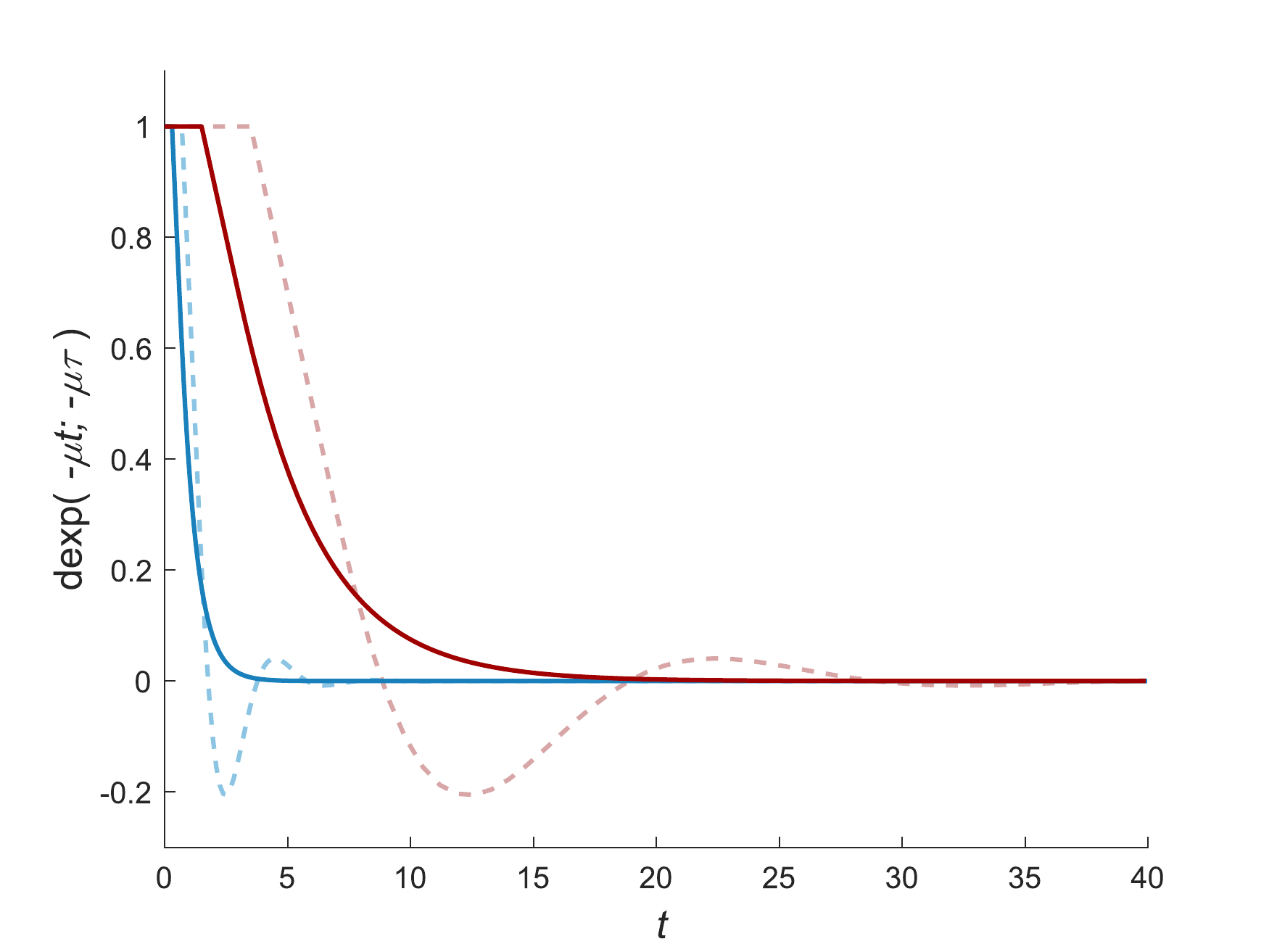}
			\caption{Examples of delay exponential functions with varying parameter settings are depicted. Blue lines represent cases where $\mu = 1$, while red lines represent $\mu = 0.2$. Solid and dashed lines correspond to different choices of $\mu\tau$: specifically, solid lines correspond to $\mu\tau = 0.3$, while dashed lines correspond to $\mu\tau = 0.7$.}
			\label{fig_dexpDistribution}
			\vspace{-1.5\baselineskip}
		\end{center}
	\end{figure}
	Examples of delay exponential functions with different choices of parameters are shown in Figure \ref{fig_dexpDistribution}. The function remains constant at $1$ for $0\leq t \leq \tau$. It then decreases linearly until $2\tau$, then for any interval $n \tau \leq t \leq (n+1)\tau$ the function will decrease at a rate proportional to $t^n$. When $\mu\tau > e^{-1}$ the function exhibits oscillatory behaviour not aligned with a survival function. Hence, we note that $\mathrm{dexp}(-\mu t;-\mu\tau)$ is a valid survival function if and only if $0 \leq \mu\tau \leq e^{-1}$. This can be seen as $\mathrm{dexp}(-\mu t;-\mu\tau)$ is right continuous, which follows directly from the definition of the Heaviside function, and $\Psi(0^+)=1$. The Laplace transform of Eq. \eqref{eq:delay_survival} is given by,
	\begin{equation}\label{eq:laplace_Phi}
		\calL_t\{\Psi(t)\}=\frac{1}{s+\mu e^{-s\tau}}.
	\end{equation}
	Using the Laplace transform and the Abelian theorem, we establish that $\lim\limits_{t\rightarrow\infty}\Psi(t)=\lim\limits_{s\rightarrow0}s\calL_t\{\Psi(t)\} = 0$. Moreover, the delay exponential function is positive and monotonically non-increasing for $0 \leq \mu\tau \leq e^{-1}$, as detailed in Appendix \ref{sec:dexp_dist_results}. These properties ensure that $\mathrm{dexp}(-\mu t;-\mu\tau)$ is a survival function. However, when $\mu\tau > e^{-1}$, delay exponential functions become negative and oscillatory, as discussed in Agarwal et al. \cite{ABBD2012}, rendering them unsuitable as probability distributions, see the dashed lines in Figure \ref{fig_dexpDistribution}.
	
	For a random variable $X$ following a delay exponential distribution, the moment generating function (MGF) is given by:
	\begin{equation}
		M_X(s)=\mathbb{E}[e^{sX}] = \frac{1}{1-\mu^{-1}se^{-s\tau}}.
	\end{equation}
	From this, we can derive its mean and variance,
	\begin{equation}
		\begin{split}
			\mathbb{E}[X] = \mu^{-1}, \quad
			\mathrm{Var}(X) = \mathbb{E}[X^2]-\mathbb{E}[X]^2 = \mu^{-2}(1-2\mu\tau).
		\end{split}
	\end{equation}
	These moments could be helpful to fit a waiting time distribution in a model. It should be noted, like in the case of the exponential distribution, the mean is inverse to the scale parameter.
	
	We will now choose our non-Markovian removal process to be a delay exponential distribution. However, we consider a waiting time of this form as it leads to delay terms in our compartment model and produces a tractable model which can be simulated exactly. Hence there may be contexts where the accuracy of model isn't as paramount as how simple it is to solve.
	
	\section{Compartment model with time delays}\label{sec:model_time_delay}
	To write a governing equation for a compartment model with a time delay, we return to Eq. \eqref{eq:gov_comp_kern} with $\Psi_i(t)$ as the delay exponential distribution. 
	Given its Laplace transform, Eq. \eqref{eq:laplace_Phi}, the memory kernel, Eq. \eqref{eq:memory_kernel}, becomes,
	\begin{equation}
		\calL_t\{ K_i(t) \}=\mu_i e^{-s\tau_i}.
	\end{equation}
	Substituting this memory kernel into the governing equation, Eq. \eqref{eq:gov_comp_kern}, we get,
	\begin{align} \label{eq:gov_comp_delay}
		\frac{d\rho_i}{dt}&=q_i^+(t)-\omega_i(t)\rho_i(t)-\Phi_i(t,0) \calL_t^{-1}\left\{\mu_i e^{-s\tau_i}\calL_t\left\{\frac{\rho_i(t)}{\Phi_i(t,0)}\right\}\right\}\\
		&=q_i^+(t)-\omega_i(t)\rho_i(t)-\mu_i\frac{\Phi_i(t,0)}{\Phi_i(t-\tau_i,0)}\rho_i(t-\tau_i)\Theta(t-\tau_i).
	\end{align}
	If in Eq. \eqref{eq:Phi_def}, $\omega_i(s)=\omega_i$, then the model simplifies to:
	\begin{equation}
		\frac{d\rho_i}{dt}=q_i^+(t)-\omega_i\rho_i(t)-\mu_i e^{-\omega_i\tau_i}\rho_i(t-\tau_i)\Theta(t-\tau_i).
	\end{equation}
	Furthermore, our initial conditions state that the compartment is empty prior to $t=0$, hence $\rho_i(t)=0$ for $t<0$. This further simplifies our governing equation for a compartment model with a delay to,
	\begin{equation}\label{eq:gov_eqn}
		\frac{d\rho_i}{dt}=q_i^+(t)-\omega_i\rho_i(t)-\mu_i e^{-\omega_i\tau_i}\rho_i(t-\tau_i).
	\end{equation}
	In this equation, the flux into the compartment, $q_i^+(t)$, can be made up of flows out of other compartments. Hence, it can include delay terms which capture the delayed removal process from another compartment. It can also include non-linear terms which may represent interactions between populations in different compartments. We now consider a method to simulate this model.
	\section{Stochastic simulation method}\label{sec:simulation}
	To simulate a compartment model with delay, we establish a stochastic simulation method for the governing equations, Eq. \eqref{eq:gov_eqn}. We use the versatile framework introduced in Xu et al. \cite{XAHBHJ2023}, which is applicable to non-Markovian compartment models. Given the non-Markovian nature of delayed removal processes, each particle in the delay compartments is assigned a waiting time sampled from the delay exponential survival function. This is achieved by numerically inverting the equation,
	\begin{equation}\label{eq:dexpSampling}
		u=\mathrm{dexp}(-\mu t;-\mu\tau),
	\end{equation}
	where $u$ is a uniform(0,1) random variable. Numerical inversion may be accomplished through standard root-finding algorithms. Efficiency can be enhanced by precomputing function values at integer multiples of $\tau$, offering a reliable initial guess. The simulation then tracks the time until the particle leaves the delay compartment as time evolves, while handling all Markovian processes similarly to the Gillespie algorithm \cite{GILLESPIE1976}. This involves considering the overall Markovian transition rates, which depend on the current system states, analogous to propensity functions in the context of chemical reactions.
	
	Consider a system of $N\ge1$ compartments with at most one delay removal process in each compartment. The dynamic state of this system can be described by the state vector $[X_1(t),\dots,X_N(t)]$, where $X_i(t)$, $1\leq i \leq N$, is the number of particles in compartment $i$ at time $t$. The exact stochastic simulation algorithm for this model is described as follows,\\
	
	\textit{Algorithm:}
	\begin{enumerate}[ {(}1{)} ]
		\item Initialize. Set the initial number of particles in each compartment. Set the system time $t=0$.
		\item{For each compartment $i$ with $n_i\geq 1$ associated Markovian removal processes, generate a waiting time $\Delta t_i$ from the survival function
			\begin{equation}\label{eq:Phi_inhomoPoisson}
				\Phi_i(\Delta t_i|t)=\exp\left(-\int_{t}^{t+\Delta t_i} \omega_i(s)X_i(t)ds\right),
			\end{equation}
			using a method of one's choice.}
		\item{If compartment $i$ includes a delayed removal process, then label each particle $k$ where $1\leq k \leq X_i(t)$. Generate a uniform(0,1) random variable, $u_1$, and then generate a waiting time $\Delta t_{i,k}$ from the delay exponential survival function, $\mathrm{dexp}(-\mu_i t;-\mu_i\tau_i)$, by numerically inverting Eq. \eqref{eq:dexpSampling}. Repeat this step for each particle $k$ in the compartment and store the list $\{\Delta t_{i,k}\}$.}
		\item{Set the time to the next transition for the system to
			\[\Delta t=\min_{\forall i,k}(\Delta t_i,\{\Delta t_{i,k}\}).\]}
		\item{When the next transition occurs in compartment $i$:
			\begin{itemize}
				\item[$\cdot$]{if $\Delta t = \Delta t_i$, use the Gillespie algorithm \cite{GILLESPIE1976} to select the $h$th Markovian removal process based on:
					\begin{equation}
						\sum_{j=1}^{h-1}\lambda_{i,j}(t)<u_2\omega_i(t)\leq\sum_{j=1}^{h}\lambda_{i,j}(t)
					\end{equation}
					where $\lambda_{i,j}(t)$ is the rate of an individual Markovian removal process, $\omega_i(t)$ is the net Markovian removal rate for compartment $i$ defined from Eq. \eqref{eq:omega_lambda}, and $u_2$ is a uniform(0,1) random variable.
				}
				\item[$\cdot$]{else if $\Delta t=\min_{\forall k}(\{\Delta t_{i,k}\})$, select the particle for removal by $\argmin_k\{\Delta t_{i,k}\}$.
				}
			\end{itemize}
			Update $X_i(t+\Delta t)$ accordingly.}
		\item{Update system time $t\leftarrow t+\Delta t$ and all waiting times: $\{\Delta t_i\}\leftarrow\{\Delta t_i - \Delta t\}$ and $\{\Delta t_{i,k}\}\leftarrow\{\Delta t_{i,k}-\Delta t\}$. Redraw $\Delta t_i$ when $\Delta t_i=0$ or if the survival function, Eq. \eqref{eq:Phi_inhomoPoisson}, has changed following the transition. Track particles entering or leaving the delayed compartment. For particles entering, generate a new waiting time using Eq. \eqref{eq:dexpSampling} and add it to $\{\Delta t_{i,k}\}$. For particles leaving, remove its waiting time from the list.}
		\item{Return to step 4 or quit.}
	\end{enumerate}
	
	We will make use of this algorithm to simulate the examples in the following section.
	
	\section{Examples}\label{sec:example}
	Compartment models with time delays have been used across many fields \cite{RHRC2018}, in this section we consider a pharmacokinetic, and an epidemiological model, and show an exact stochastic simulation for each model. We use the general framework established in this paper to write the governing equations of the specific examples.
	The pharmacokinetic model we consider a drug absorption, with a delayed clearance term. For the epidemiological model we consider is an SIS model with a delayed re-susceptibility term.

	\subsection{Pharmacokinetic model with clearance delay}\label{sec:pharma}
	Pharmacokinetic models describe the absorption, clearance and effects of particular drugs in a system. When drugs are administered, there is often a delay observed between when the drug is administered and when the drug is absorbed. This observed delay has lead to the incorporation of delay terms into pharmacokinetic models \cite{RHRC2018}. A delay in a pharmacokinetic system can also occur due to a delayed drug clearance process. 
	
	Here we present a simple two compartment model obtained from Eq. \eqref{eq:gov_eqn}. A drug is released into a transport compartment $\rho_1(t)=x(t)$, before making its way to $\rho_2(t)=A(t)$, the target location for the drug. From $A(t)$ the drug is cleared from the system following a delayed process. In this model, $x(t)$ and $A(t)$ represent the concentration of the drug in each compartment. Starting with the dynamics of $x(t)$, there is a Markovian clearance out of the compartment, so $\omega_1=k$. There are no other flows in or out of $x(t)$, hence $\tau_1=0$, $\mu_1=0$ and $q^+_1(t)=0$. So, the governing equation for $x(t)$ is,
	\begin{equation}
		\frac{dx}{dt}=-kx(t).
	\end{equation}
	The flow into the target compartment $A(t)$ is the outflow of $x(t)$, hence $q^+_2(t)=kx(t)$. There is no Markovian clearance from $A(t)$ so $w_2=0$. The clearance from $A(t)$ has a $\tau_2=\tau$ delay. The systematic clearance rate is $C$ and the volume of $A(t)$ is $V$ which we have taken to be constant. Together they are used to define the scale parameter, $\mu_2=C/V$. Hence the governing equation for $A(t)$ is,
	\begin{equation}
		\frac{dA}{dt}=kx(t)-\frac{C}{V}A(t-\tau).
	\end{equation}
	Both compartments are empty for $t<0$, then a  dose of the drug, $x_0$, is injected into the transport compartment at $t=0$. So $x(t)=x_0\delta(t)$ and $A(t)=0$ for $t\in [-\tau,0]$ form our initial conditions for the model.
	
	While our model has only one transport compartment, this model can be generalized to include $n$ transport compartments.
	
	\begin{figure}[!htbp]
		\begin{center}
			\begin{subfigure}{.75\textwidth} 
				\includegraphics[width=1\linewidth]{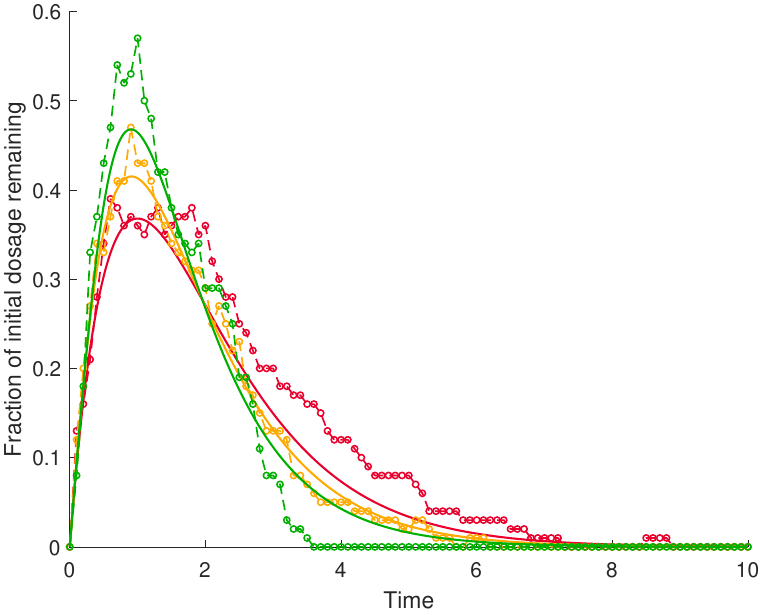}
				\caption{}	
				\label{fig:PK_sub1}
			\end{subfigure}
			\begin{subfigure}{.75\textwidth}
				\centering
				\includegraphics[width=1\linewidth]{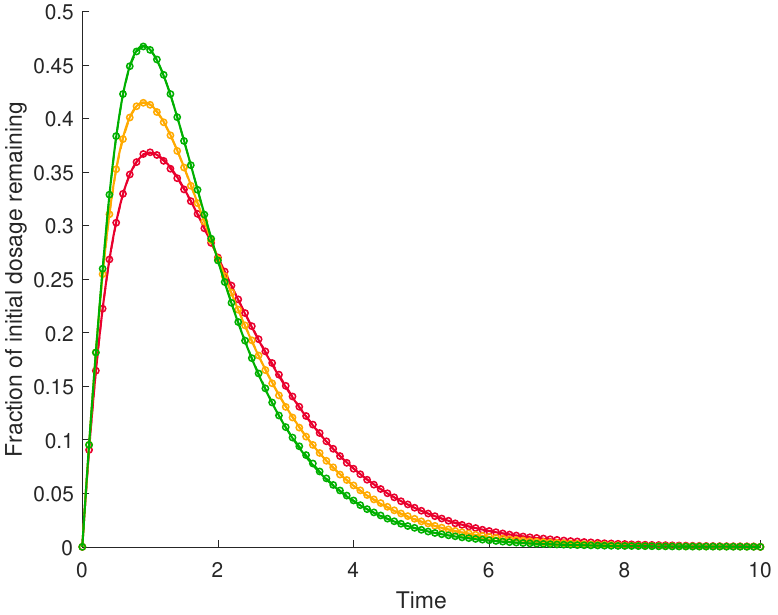}
				\caption{}
				\label{fig:PK_sub2}
			\end{subfigure}
			\caption{Stochastic and deterministic solutions of the two-compartment pharmacokinetic model with a delayed clearance. We consider a drug transport rate $k=1$, a clearance rate $C/V = 1$ and three time delays $\tau = 0$ (red), $0.2$ (orange) and $0.35$ (green), with an initial injection $x_0=100$. (a) Representative sample paths of the stochastic simulation (circle--dashed lines) and the corresponding deterministic solutions (solid lines) are shown. (b) $2000$ sample paths are averaged to show the convergence between the stochastic and the deterministic solution.} 
			\label{fig_PKs}
			\vspace{-1.5\baselineskip}
		\end{center}
	\end{figure}
	
	We simulate the model with parameters $k=1$ and $C/V=1$, starting with an initial injection of $x_0=100$. We consider the delay, $\tau$ at $0$, $0.2$ and $0.35$. In Figure \ref{fig_PKs} we show the concentration of the drug in the target compartment and consider the impact of different delay parameters. The simulation reveals that a larger delay in the drug clearance process results in a higher peak drug concentration in the target compartment. However larger clearance delays also lead to the target compartment clearing faster.
	
	In Figure \ref{fig:PK_sub1}, we display sample paths of the stochastic simulations alongside the corresponding deterministic solution. In Figure \ref{fig:PK_sub2}, we demonstrate the convergence between the stochastic and deterministic solutions, by averaging $2000$ sample paths. 
	
	We note that the stochastic solution converges to the deterministic solution even with a small population size. This will be the case for our simulation algorithm when the model contains only linear terms.

	\subsection{SIS model with delayed re-susceptibility}\label{sec:SIS}
	The SIS epidemic model with a delayed re-susceptibility is a generalisation of the standard SIS model and is a special case of Hethcote and van den Driessche \cite{HD1995}. In this model the population is split into two compartments: the Susceptible ($S$) and Infected populations ($I$). Individuals are born into the Susceptible compartment and they are infected through interaction with infected individuals. Once they recover from the disease, they move back to the Susceptible compartment and can be re-infected. Individuals can leave the system from either compartment through death. This model is obtained from Eq. \eqref{eq:gov_eqn} by setting $\rho_1(t)=S(t)$ and $\rho_2(t)=I(t)$. The flux into the infected compartment is, $q^+_2(t)=\beta S(t)I(t)$. We consider a process where the whole population shares the same death rate. So the Markovian removal process of the infected population is, $\omega_2=d$, where $d>0$. The re-susceptibility of the population will be a delayed process hence, $\mu_2=\gamma$ and $\tau_2=\tau$. This gives us the governing equation for the infected compartment,
	\begin{equation}\label{eq:dSIS_I}
		\frac{dI}{dt}=\lambda S(t)I(t)-\gamma e^{-d\tau} I(t-\tau)-dI(t).
	\end{equation}
	To write the governing equation for the susceptible compartment, we take the flux into the compartment to be the individuals who become re-susceptible, as well as the births of the population. The birth term is proportional to the size of the whole population. Hence, $q^+_1(t)=b(S(t)+I(t))+\gamma e^{-d\tau} I(t-\tau)$, where $b>0$ is the birth rate in the population. Individuals leave the susceptible population either through death or infection. Hence $\omega_1=d+\beta I(t)$. There is no delayed removal process out of the susceptible compartment, hence $\tau_1=0$ and $\mu_1=0$. This yields the governing equation for the susceptible compartment,
	\begin{equation}\label{eq:dSIS_S}
		\frac{dS}{dt}=b(S(t)+I(t))+\gamma e^{-d\tau} I(t-\tau)-\lambda S(t)I(t)-dS(t).
	\end{equation}
	To complete the dynamics of the model we set $S(0)=s_0$ and $I(t)=i_0\delta(t)$, for $t\in [-\tau,0]$.
	\begin{figure}[!htbp]
		\begin{center}
			\begin{subfigure}{.75\textwidth} 	
				\includegraphics[width=1\linewidth]{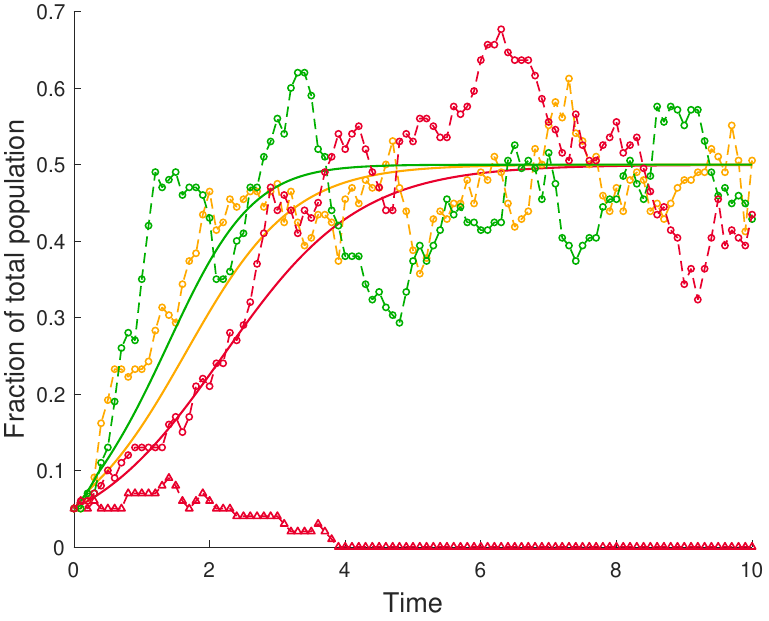}
				\caption{}
				\label{fig:SIS_sub1}
			\end{subfigure}
			\begin{subfigure}{.75\textwidth}
				\centering
				\includegraphics[width=1\linewidth]{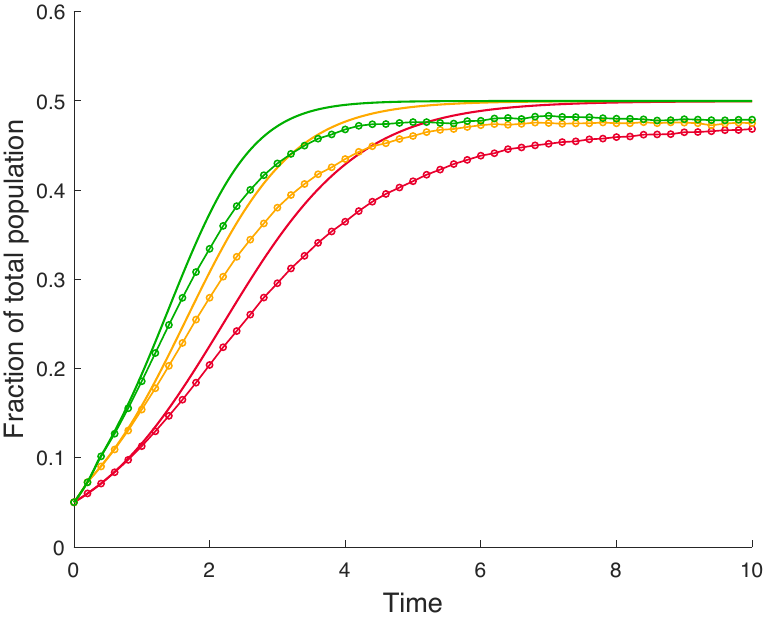}
				\caption{}
				\label{fig:SIS_sub2}
			\end{subfigure}
			\caption{Stochastic and deterministic solutions of the SIS model with delayed re-susceptibility. We consider the re-susceptibility process with rate $\gamma = 1$ three time delays $\tau = 0$ (red), $0.2$ (orange) and $0.35$ (green), for a small population with initial populations $s_0=95$ and $i_0=5$. (a) Representative sample paths of the infected population (circle--dashed lines) and the corresponding deterministic solutions (solid lines) are presented as fractions of the total population $P(t) = S(t)+I(t)$. The triangle-dashed line represents a sample path leading to disease extinction. (b) The deterministic solution is shown against a stochastic solution, calculated as an average of $2000$ sample paths.} 
			\label{fig_SISs}
			\vspace{-1.5\baselineskip}
		\end{center}
	\end{figure}
	
	\begin{figure}[!htbp]
		\begin{center}
			\begin{subfigure}{.75\textwidth} 	
				\includegraphics[width=1\linewidth]{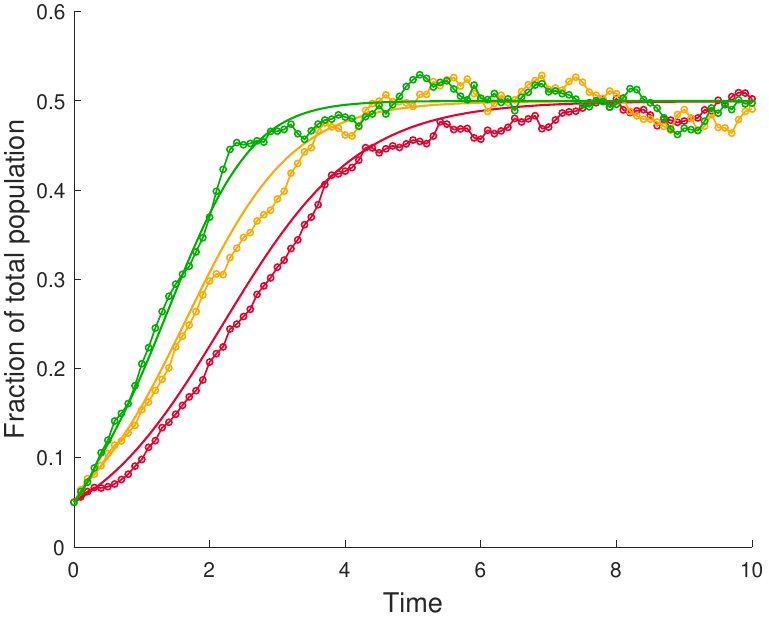}%to be updated
				\caption{}
				\label{fig:SIS_sub1_large}
			\end{subfigure}
			\begin{subfigure}{.75\textwidth}
				\centering
				\includegraphics[width=1\linewidth]{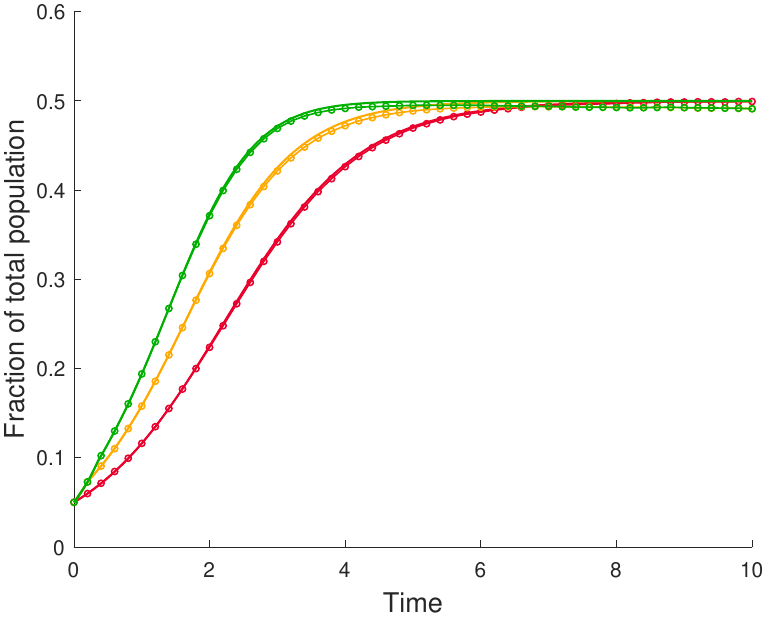}
				\caption{}
				\label{fig:SIS_sub2_large}
			\end{subfigure}
			\caption{Stochastic and deterministic solutions of the SIS model with delayed re-susceptibility. We consider the re-susceptibility process with rate $\gamma = 1$ three time delays $\tau = 0$ (red), $0.2$ (orange) and $0.35$ (green), for a large population with initial populations $s_0=1900$ and $i_0=100$. (a) Representative sample paths of the infected population (circle--dashed lines) and the corresponding deterministic solutions (solid lines) are presented as fractions of the total population $P(t) = S(t)+I(t)$. (b) The deterministic solution is shown against a stochastic solution, calculated as an average of $2000$ sample paths.} 
			\label{fig_SISs_large}
			\vspace{-1.5\baselineskip}
		\end{center}
	\end{figure}
	
	To simulate this model we have set the parameters to $b=0.1$, $d=0.1$, and $\lambda = 2/P_0$, where $P_0 =s_0+i_0$ denotes the total initial population. We will consider two choices for $P_0$ and take the delay $\tau$ at $0$, $0.2$ and $0.35$.  In Figure \ref{fig_SISs} we consider a small initial population where $P_0=100$. We show the population of the Infective compartment, $I$ over time. In Figure \ref{fig:SIS_sub1} we consider sample paths for each delay shown bye the circle-dashed lines. We note that the non-linear flux term, $q^+_2(t)=\beta S(t)I(t)$, leads to a high probability of extinction events, as shown by the triangle-dashed line in the figure. In Figure \ref{fig:SIS_sub2} we show the deterministic and stochastic solutions. The stochastic solution is taken as the average of $2000$ paths. We observe a significant deviation between the stochastic and deterministic solutions with a small initial population, this is another result of the non-linear flux term.  From the simulation, we observe that a delayed re-susceptibility affects the time taken for the system to reach the endemic states. Specifically, the larger the delay, the faster the system approaches the endemic states. However, the delay does not alter the endemic states themselves as seen from the deterministic solution shown by the solid lines.
	
	We show much better agreement between the stochastic and deterministic solutions with a larger initial population, $P_0=2000$, in Figure \ref{fig_SISs_large}. We have kept the other parameters consistent with the ones used in Figure \ref{fig_SISs}. We show sample paths for the larger population in Figure \ref{fig:SIS_sub1_large}. While the probability of extinction events is not zero, it is much less likely with a larger population. In this plot we note that the sample paths are much more consistent with the deterministic solution. By averaging $2000$ sample paths for each delay, we observe the convergence of the deterministic and stochastic solutions in Figure \ref{fig:SIS_sub2_large}.

	\section{Discussion}\label{sec:summary}
	Compartment models with delays have been observed across a wide range of fields. Typically they are derived in an ad hoc fashion, justified by considering an incubation or delayed process. Here, we have derived the governing equations for a compartment model with delayed terms from an underlying stochastic process. The resulting dynamics are represented by the master equation, Eq. \eqref{eq:gov_eqn}. We first derived a master equation for a compartment model with a general non-Markovian removal process through Section \ref{sec:deriv}.
	
	We then showed that a delay exponential function is a probability distribution under particular constraints in Section \ref{sec:dexp_dist}. We recovered a compartment model with delayed removal terms when we took the delay exponential distribution as the waiting time in the master equations in Section \ref{sec:model_time_delay}.
	
	We provided a simulation method in Section \ref{sec:simulation} for the compartment model with delays. We establish an exact stochastic simulation method that can solve DDEs and provides stochastic simulations to a range of compartment models. The method is efficient and highly tractable. 
	
	Finally, we provided some examples of compartment models with delays in section \ref{sec:example}. For these examples the governing evolution equations were obtained and the simulations were done using the methods set out in this paper. We numerically demonstrated that, the deterministic solution aligns with stochastic simulation results with a small population when all flux terms are linearly dependent on the state of the current compartment or deterministic. However, for small populations, divergence may occur when there is a non-linear flux term, as observed in the SIS model with delayed resusceptibility. We showed convergence in the limit with large populations when we had a non-linear flux. The convergence between deterministic and stochastic solutions has been studied by Kurts et al. for Markov processes \cite{KurtzBook1986}, and by Armbruster et al. within the context of epidemiological models \cite{Armbruster2017a,Armbruster2017b}, similar studies on delay models are yet to be investigated.
	
	In this paper we considered the underlying stochastic process that leads to a delay term in a compartment model. This provides greater insight into the mechanisms at play when compartment models with delays are chosen. However our model does not account for all occurrences of delay terms in compartment models.
	
	We have focused on delays that come from the flux out of compartments. Our stochastic process does not capture the dynamics of every existing DDE compartment models. The DDEs we derive are computationally tractable and amenable to exact stochastic simulations. The physical interpretation of parameters in our models are clear. Our formulation will allow for further generalizations of DDEs in models via the modification of the underlying stochastic process.
		\section*{Acknowledgements}
	This work was supported by the Australian Research Council grant number DP200100345.
	
	We give a special thank you to productive conversations with Bruce Henry, Daniel Han and Boris Huang.
	\appendix
	\section{On positiveness and complete monotonicity of the delay exponential survival function}\label{sec:dexp_dist_results}
	Consider the delay differential equation,
	\begin{equation}
		\frac{d\rho(t)}{dt} = - \mu\rho(t-\tau), \quad \mu,t, \tau \in \mathbb{R}^+,
	\end{equation}
	with initial condition,
	\begin{equation}
		\rho(t) = \phi(t) \quad\text{for} \quad t\in [-\tau, 0]. 
	\end{equation}
	We know that the fundamental solution with initial condition $\phi(t)=0$ for $t<0$ and $\phi(0)=1$ is given by,
	\begin{equation}
		\rho(t) = \text{dexp}(-\mu t; -\mu\tau).
	\end{equation}
	We also know that when $\mu\tau<\frac{1}{e}$, the characteristic function, 
	\begin{equation}
		\lambda + \mu e^{-\lambda\tau}=0,
	\end{equation}
	has negative roots given by $\frac{W_0(-\mu\tau)}{\tau}$ and $\frac{W_{-1}(-\mu\tau)}{\tau}$, where $W_k(z)$ is the $k$th branch of the Lambert W function \cite{CGHJK1996} defined by
	\begin{equation}
		W_k(z)e^{W_k(z)}=z,
	\end{equation}
	with
	\begin{equation}
		\lim_{|z|\rightarrow\infty}W_k(z)\rightarrow \mathrm{ln}_kz = \mathrm{ln}z + 2\pi ik.
	\end{equation}
	Therefore, there exist positive completely monotonic solutions, for example,
	\begin{equation}
		u(t) = e^{\frac{W_0(-\mu\tau)}{\tau}t}, \quad t\in[-\tau,\infty).
	\end{equation}
	Suppose that $\text{dexp}(-\mu t; -\mu\tau)$ is not strictly positive, and assume that the earliest time it becomes zero is at $t=t_0$, vis., $\text{dexp}(-\mu t_0; -\mu\tau)=0$. It is obvious that $t_0>\tau$. The particular solution $u(t)$ is related to the fundamental solution via:
	\begin{equation}
		u(t) = \phi(-\tau)\text{dexp}(-\mu (t+\tau); -\mu\tau)+\int_{0}^{\tau}\phi'(s-\tau)\text{dexp}(-\mu (t+\tau-s); -\mu\tau)ds,
	\end{equation}
	where $\phi'(t) = \frac{W_0(-\mu\tau)}{\tau}e^{\frac{W_0(-\mu\tau)}{\tau}t}$. Then we see that the value of $u(t)$ at $t = t_0-\tau$ is given by,
	\begin{equation}\label{eq_utodexp}
		\begin{split}
			u(t_0-\tau) &= \phi(-\tau)\text{dexp}(-\mu t_0; -\mu\tau)+\int_{0}^{\tau}\phi'(s-\tau)\text{dexp}(-\mu (t_0-s); -\mu\tau)ds,\\
			&=\int_{0}^{\tau}\phi'(s-\tau)\text{dexp}(-\mu (t_0-s); -\mu\tau)ds.
		\end{split}
	\end{equation}
	Since $\phi'(t)<0$ for $t\in[-\tau,0]$, and $\text{dexp}(-\mu t; -\mu\tau)\geq0$ for $t\in[t_0-\tau,t_0]$, the integral on the RHS of \eqref{eq_utodexp} is less or equal to zero. This contradicts the condition of $u(t)$ being a strictly positive solution for all $t$.
	
	The above shows that $\text{dexp}(-\mu t; -\mu\tau)$ is strictly positive for $t\geq-\tau$. The complete monotonicity then follows from the property:
	\begin{equation}
		\frac{d}{dt}\text{dexp}(-\mu t; -\mu\tau)= -\mu\text{dexp}(-\mu (t-\tau); -\mu\tau),
	\end{equation}
	and hence,
	\begin{equation}
		(-1)^m\frac{d^m}{dt^m}\text{dexp}(-\mu t; -\mu\tau)= \mu^m\text{dexp}(-\mu (t-m\tau); -\mu\tau)\geq0.
	\end{equation}
	The above proof works explicitly for $\mu\tau<1/e$. For $\mu\tau>1/e$, $\text{dexp}(-\mu t; -\mu\tau)$ becomes oscillatory as demonstrated in Agarwal et al. \cite{ABBD2012}. This shows that the condition for completely monotonic $\text{dexp}(-\mu t; -\mu\tau)$ is given by $\mu\tau<1/e$.

	\vfill
	\bibliographystyle{unsrt}  
	\bibliography{delay_compartment_models}  
	
\end{document}